\documentclass[12pt,a4paper]{article}

\usepackage{cmap}
\usepackage[utf8x]{inputenc}
\usepackage[T1]{fontenc}
\usepackage[english]{babel}
\usepackage{amssymb,amsmath}
\usepackage{amsthm}
\usepackage{bbm}
\usepackage{microtype}
\usepackage{lmodern}

{
\rmfamily
\DeclareFontShape{T1}{lmr}{m}{sc}{<->ssub*cmr/m/sc}{}
\DeclareFontShape{T1}{lmr}{b}{sc}{<->ssub*cmr/b/sc}{}
\DeclareFontShape{T1}{lmr}{bx}{sc}{<->ssub*cmr/bx/sc}{}
}

\newcommand{\thmheadercommand}[1]{\textbf{\scshape{}#1.\\*}}

\newtheoremstyle{yannthm}{\topsep}{\topsep}{\slshape}{}{\scshape\bfseries}{.}{.5em}{%
\thmname{#1}\thmnumber{ #2}\thmnote{#3}%
}
\newtheoremstyle{yannthm2}{\topsep}{\topsep}{}{}{\scshape\bfseries}{.}{.5em}{%
\thmname{#1}\thmnumber{ #2}\thmnote{#3}%
}

\def\d{\operatorname{d}\!{}}

\def\Z{{\mathbb{Z}}}

\renewcommand{\geq}{\geqslant}
\renewcommand{\leq}{\leqslant}

\renewcommand{\epsilon}{\varepsilon}
\renewcommand{\phi}{\varphi}

\DeclareMathOperator{\diam}{diam}

\let\oldPr\Pr
\renewcommand{\Pr}{\oldPr\nolimits}

\newcommand{\abs}[1]{\left\lvert#1\right\rvert}
\newcommand{\norm}[1]{\left\lVert#1\right\rVert}

\newcommand{\presgroup}[2]{\left\langle\,#1 \mid  #2\,\right\rangle}

\newenvironment{dem}[1][]{\begin{proof}[\thmheadercommand{Proof#1}]~\newline\ignorespaces}{\end{proof}}

{
\theoremstyle{yannthm}
\newtheorem{defi}{Definition}
\newtheorem*{defi*}{Definition}

\newtheorem*{prop*}{Proposition}
\newtheorem{thm}[defi]{Theorem}
\newtheorem*{thm*}{Theorem}
\newtheorem{lem}[defi]{Lemma}
\newtheorem*{lem*}{Lemma}
\newtheorem{cor}[defi]{Corollary}
\newtheorem*{cor*}{Corollary}

\newtheorem*{ex*}{Example}

\newtheorem*{subenonce}{}

\theoremstyle{yannthm2}

\newtheorem*{exo*}{Exercise}

\newtheorem{rem}[defi]{Remark}
\newtheorem*{rem*}{Remark}

\newtheorem*{subenonce2}{}

}

\usepackage{graphicx}

\def\d{\partial}

\title {On a small cancellation theorem of Gromov}

\author{Yann Ollivier}

\date{}

\newcommand{\tempgras}[1]{#1}

\begin{document}

\maketitle

\begin{abstract}
We give a combinatorial proof of a theorem of Gromov, which extends the
scope of small cancellation theory to group presentations arising from
labelled graphs.
\end{abstract}

In this paper we present a combinatorial proof of a small cancellation
theorem stated by M.~Gromov in~\cite{G3}, which strongly generalizes
the
usual tool of small cancellation. Our aim is to complete the
six-line-long proof given in~\cite{G3} (which invokes geometric
arguments).

Small cancellation theory  is an easy-to-apply tool of combinatorial
group theory (see~\cite{Sch} for an old but nicely written introduction,
or~\cite{GH} and~\cite{LS}). In one of its forms, it basically asserts
that if we face a group presentation in which no two relators share a
common subword of length greater than $1/6$ of their length, then the
group so defined is hyperbolic (in the sense of~\cite{Gro1},
see also~\cite{GH} or~\cite{S} for basic properties), and infinite
except for some trivial cases.

The theorem extends these conclusions to much more general situations.
Suppose that we are given a finite graph whose edges are labelled by
generators of the free group $F_m$ and their inverses (in a reduced way,
see definition below). If no word of length greater than $1/6$
times the length of the smallest loop of the graph appears twice on the
graph, then the presentation obtained by taking as relations all the
words read on all loops of the graph defines a hyperbolic group which
(if the rank of the graph is at least $m+1$, to avoid trivial cases) is
infinite.  Moreover, the given graph naturally embeds isometrically into
the Cayley graph of the group.

The new theorem reduces to the classical one when the graph is a disjoint
union of circles.  Noticeably, this criterion is as easy to use as the
standard one.

For example, ordinary small cancellation theory cannot deal with such
simple group presentations as $\presgroup{S}{w_1=w_2=w_3}$ because the
two relators involved here, $w_1w_2^{-1}$ and $w_1w_3^{-1}$, share a long
common subword. The new theorem can handle such situations: for
``arbitrary enough'' words $w_1,w_2,w_3$, such presentations will define
infinite, hyperbolic groups, although from the classical point of view
these presentations satisfy (e.g.\ if the $w_i$'s have the same length) a
priori only the $C'(1/2)$ condition from which nothing could be deduced.

\bigskip

The groups obtained by this process can in some cases be noticeably
different from ordinary small cancellation groups. For example, the
graphs used by Gromov in~\cite{G3} provide groups having Kazhdan's
property $(T)$ (see~\cite{Sil03}), whereas ordinary small cancellation
groups cannot have property $(T)$ (see~\cite{Wis}).

Most importantly, this technique allows to (quasi-)embed prescribed
gra\-phs into the Cayley graphs of hyperbolic groups. It is the basic
construction involved in the announcement of a counter-example to the
Baum-Connes conjecture with coefficients (see~\cite{HLS} which elaborates
on~\cite{G3}, or~\cite{Gh} for a survey).  Indeed, this counter-example
is obtained by constructing a finitely generated group (which is a limit
of hyperbolic groups) whose Cayley graph quasi-isometrically contains an
infinite family of expanders.

Moreover, this technique will be used in~\cite{OW} to construct new
examples of groups with property~$(T)$.

\section{Statement and discussion}

Let $S$ be a finite set, in which an involution without fixed point,
called \emph{being inverse}, is given.  The elements of $S$ are called
\emph{letters}.

A \emph{word} is a finite sequence of letters. The inverse of a
word is the word made of the inverse letters put in reverse order.
A word is called \emph{reduced} if it does not contain a letter
immediately followed by its inverse.

A \emph{labelled graph} is an unoriented graph in which each unoriented
edge is considered as a couple of two oriented edges, and each oriented
edge bears a letter such that opposite edges bear inverse letters. We
require maps of labelled graphs to preserve the labels.

A labelled graph is said to be \emph{reduced} if there is no pair of
oriented edges arising from the same vertex and bearing the same letter.

Note that a word can be seen as a (linear) labelled graph, which we will
implicitly do from now on. The word is reduced if and only if the labelled
graph is.

A \emph{piece} of a labelled graph is a word which has two different
immersions in the labelled graph. (An immersion is a locally injective
map of labelled graphs. Two immersions are considered different if
they are different as maps.) This is analogue to the
traditional piece of small cancellation theory.

A \emph{standard family of cycles} for a connected graph is a set of
paths in the graph, generating the fundamental group, such that there
exists a maximal subtree of the graph such that, when the subtree is
contracted to a point (so that the graph becomes a bouquet of circles),
the set of generating cycles is exactly the set of these circles. There
always exists some. If the graph is not connected, a standard family of
cycles is one which is standard on each component.

A \emph{generating family of cycles} is a family of cycles generating the
fundamental group of each connected component of the graph
(maybe up to adding initial and final segments joining these cycles to
some basepoint).

A graph is \emph{non-filamenteous} if every edge belongs to some immersed
cycle.

\bigskip

We are now in a position to state the theorem.

\begin{thm}[ (M.~Gromov, \cite{G3})]
\label{thm:main}
Let $\Gamma$ be a finite reduced non-filamenteous labelled graph. Let $R$ be the set of words read on
all cycles of $\Gamma$ (or on a generating family of cycles). Let $g$ be the girth of $\Gamma$ and $\Lambda$ be the
length of the longest piece of $\Gamma$.

If $\Lambda<g/6$ then the presentation $\presgroup{S}{R}$ defines a
group $G$ enjoying the following properties.
\begin{enumerate}
\item It is hyperbolic, torsion-free.
\item Any presentation of $G$ by the words read on a standard family of
cycles of $\Gamma$ is aspherical (in the sense of
Definition~\ref{def:asphericity}), hence the cohomological dimension
of $G$ is at most $2$.
\item The Euler characteristic of $G$ is $\chi(G)=1-\abs{S}/2+b_1(\Gamma)$.
In particular, if the rank of the fundamental group of $\Gamma$ is greater than
the number of generators, $G$ is infinite and not quasi-isometric to
$\mathbb{Z}$.
\item The shortest relation in $G$ is of length $g$.
\item For any reduced word $w$ representing the identity in $G$, some cyclic
permutation of $w$ contains a subword of a word read on a circle immersed
in $\Gamma$, of length
at least $(1-3\Lambda/g)$ (which is more than $1/2$) times the length of this cycle.
\item The natural maps from each connected component of the labelled
graph $\Gamma$ into the Cayley graph of $G$ are isometric embeddings.
\end{enumerate}
\end{thm}

If $\Gamma$ is a disjoint union of circles, this theorem almost reduces
to ordinary $1/6$ small cancellation theory. The ``almost'' accounts for
the fact that the length of a shared piece between two relators is
supposed to be less than $1/6$ the length of the smallest of the two
relators in ordinary small cancellation theory, and less than $1/6$ the
length of the smallest of all relators in our case; this is handled
through the following remark (which we do not prove in order not to have
still heavier notation).

\begin{rem} It is clear from
the proof that the assumption in the theorem can be replaced by the
following slightly weaker one: for each piece, its length is less than
$1/6$ times the length of any cycle of the graph on which the piece
appears.
\end{rem}

With this latter assumption, the theorem reduces to ordinary small cancellation
when the graph is a disjoint union of circles.

\begin{rem}
Non-filamenteousness is needed only to ensure isometric embedding of the
graph (filaments may not embed isometrically if $\Lambda\geq g/8$).
\end{rem}

The group obtained is not always non-elementary: for example, if there
are three generators $a,b,c$ and the graph consists in two points joined
by three edges bearing $a$, $b$ and $c$ respectively, one obtains the
presentation $\presgroup{a,b,c}{a=b=c}$ which defines $\Z$.
However, since the cohomological dimension is at most $2$, it is easy to
check (computing the Euler characteristic) that if the rank of the
fundamental group of $\Gamma$ is greater than the number of generators, then
$G$ is non-elementary.

\medskip

This theorem is not stated explicitly in~\cite{G3} in the form we give
but using a much more abstract and more powerful formalism of ``rotation
families of groups'' (\cite{G3}, section~2).  In the vocabulary thereof,
the case presented here is when this rotation family contains only one
subgroup of the free group (and its conjugates), namely the one generated
by the words read on cycles of the graph with some base point; the
corresponding ``invariant line'' $U$ is the universal cover of the
labelled graph $\Gamma$ (viewed embedded in the Cayley graph of the free
group). Reducedness of the labelling ensures convexity.

Elements for a proof of the theorem for very small values of $\Lambda/g$
(instead of $\Lambda/g<1/6$) using geometric rather than combinatorial
tools, can be found in~\cite{G2} (see also~\cite{G3}, p.~88).

In~\cite{G3}, this theorem is applied to a random labelling (or rather
a variant, Theorem~\ref{thm:variantmain} below, in which reducedness is
replaced with quasi-geodesi\-city). It is not difficult, using for
example the techniques described in~\cite{O}, to check that a random
labelling satisfies the small cancellation and quasi-geodesi\-city
assumptions.

\paragraph{Acknowledgements.} I would like to
thank Thomas Delzant for having brought the problem to my attention
and for very careful reading, Étienne Ghys, Pierre Pansu and Alain
Valette for helpful discussions and comments on the manuscript, and
Mikhail Ostrovskii for pointing out an error in an earlier version of the
proof of the isometric embedding property.

\section{Idea of proof}

The line of the argument is as follows: Choose a presentation of $G$ by
the words read on a standard generating family of cycles of $\Gamma$. We
will study the isoperimetry of van Kampen diagrams with respect to this
set of relations: we will show that the number of faces in such
diagrams is linearly bounded by its boundary length.

Define a labelled complex $\Gamma_2$ by attaching to $\Gamma$ a disk for each
cycle in the family. Now each face of a van Kampen diagram for this presentation
can be lifted (in a unique way) to $\Gamma_2$. For any edge between two
faces of the diagram, either these two faces are already adjacent along ``the
same'' edge in $\Gamma_2$ or they are not.

Decompose the diagram into maximal parts all edges of which originate
from $\Gamma_2$ in this sense. Now gluings between these parts do not
originate from $\Gamma_2$ and thus constitute pieces. So these parts are
in classical $1/6$ small cancellation with respect to each other, and so
the boundary length of the diagram is controlled in terms of the boundary
lengths of these parts. We get the other usual consequences of small
cancellation theory as well (asphericity, radius of injectivity...).
Technicalities arise from the necessity to perform some so-called
``diamond moves'' and from the maybe non-simple connectedness of these
parts.

To reach the conclusion it is then enough to work inside each part. Since
each part lifts to $\Gamma_2$ its boundary word is the word read on some
null-homotopic cycle in $\Gamma_2$. So this cycle is the product of
elements our generating family of cycles, and for isoperimetry we have to
control the number of terms in this product (the number of faces in the
part) in function of the length of the cycle (the boundary length of the
part). This is achieved by decomposing the considered cycle into a
product of cycles shorter than three times the diameter of the graph. As
there are only finitely many such short cycles we are done.

\section{Proof (expanded version)}

We now give some more definitions which are useful for the proof.

\begin{defi}
A \emph{labelled complex} is a finite unoriented combinatorial
$2$-complex the interior of every face of which is homeomorphic to an
open disk $D_{n+1}$ with $n\geq 0$ holes ($n$ depends on the face), 
such that its $1$-skeleton is equipped with a labelled graph structure.

A labelled complex is said to be \emph{reduced} if its $1$-skeleton is.
\end{defi}

Each face of such a complex defines a set of
\emph{\tempgras{contour} words}: If the interior of the face is
homeomorphic to an open disk $D_{n+1}$ with $n$ holes, the contour words
are the $n+1$ cyclic words read by moving around the $n+1$
\tempgras{boundary} components of $D_{n+1}$. The words in this set are
considered as oriented cyclic words, and counted with multiplicities.

We require a map of labelled complexes to preserve labels (but it may
change orientation of faces, sending a face to a face with inverse
\tempgras{contour} labels --- this amounts to considering maps between the
corresponding oriented complexes).

\begin{defi}
A \emph{tile} is a planar labelled complex with only one
face (not necessarily simply connected) and each edge of which belongs to
the combinatorial boundary of the face with multiplicity one.  We do not fix the
embedding in the plane.
\end{defi}

It follows from the definition that the \tempgras{contour} of a tile
coincides with its \tempgras{boundary}.

By our definition of maps between labelled complexes, a tile is
considered equal to the tile bearing the inverse \tempgras{boundary}
words. 

Convention: A tile may bear a word which is not simple (i.e.\ is a power
of a smaller word). In this case the tile would have a non-trivial
automorphism. To prevent this, say that on each \tempgras{boundary}
component of a tile we mark a starting point and that a map between tiles
has to preserve marked points. This is useful for the study of
asphericity and torsion (see Definition~\ref{def:asphericity}).

To any planar labelled complex with only one face we can associate a tile
in the following way: First, remove the edges that do not belong to the
adherence of the interior of the complex (the ``filaments''). Then, the
obtained one-face complex immersed in the plane is the image of some
one-face complex embedded in the plane by a cellular map (this complex is
constructed by ungluing along the internal edges). This is an embedding in
the plane of some tile, which we call the \emph{tile associated to} the
one-face labelled complex.

\begin{defi}
\label{def:tileofcomplex}
A \emph{tile of a labelled complex} is the tile associated to any of its
faces.
\end{defi}

The \emph{length} of a tile is the length of its \tempgras{boundary}.

\begin{defi}
A \emph{piece} with respect to a set of tiles is a word which has
immersions in the \tempgras{boundary} of two different tiles, or two distinct
immersions in the \tempgras{boundary} of one tile.
\end{defi}

\begin{defi}
A \emph{puzzle} with respect to a set of tiles is a planar
labelled complex all tiles of which belong to this set of
tiles (the same tile may appear several times in a puzzle). The \emph{set
of boundary words} of a puzzle is the set of words read on its
boundary components (with multiplicities and orientations).

A \emph{spherical puzzle} is the same drawn on a sphere instead of the
plane, that is, a labelled complex which is a combinatorial $2$-sphere,
all tiles of which belong to this set of tiles.

A puzzle is said to be \emph{minimal} if it has the minimal number of
tiles among all puzzles having the same set of \tempgras{boundary} words.

A puzzle is said to be \emph{van Kampen-reduced} if there is no pair of
adjacent faces such that the words read on the external \tempgras{contour} of these two
faces are inverse and the position (with respect to the marked point) of
the letter read at a common edge of these faces is the same in the two
copies of the \tempgras{contour} word of these faces.
\end{defi}

So a puzzle is roughly speaking a van Kampen diagram in which we allow
non-simply connected faces. The last definition given corresponds to
reduced van Kampen diagrams (see~\cite{LS}).  (Incidentally, a reduced
puzzle is van Kampen-reduced, though the converse is not necessarily
true.)

\begin{defi}
\label{def:asphericity}
A presentation of a group is said to be \emph{aspherical} if the set of
tiles whose \tempgras{boundary} words are the relators of the
presentation admits no van Kampen-reduced spherical puzzle.
\end{defi}

There are several notions of aspherical presentations in the literature
(see e.g.~\cite{CCH} for five of them).  Our definition of asphericity
coincides with the one in~\cite{Gersten87}, p.~31 (in which asphericity
is termed ``every spherical diagram is diagrammatically reducible'').
It is thus stronger than the
one(s) in~\cite{LS}, the main difference being that we mark a starting
point on the \tempgras{boundary} of each tile (see the discussion
in~\cite{Gersten87}). In particular, with our (and \cite{Gersten87}'s,
contra~\cite{LS}) convention, a presentation such as
$\presgroup{S}{w^n=1}$ (with $n\geq 2$) is not aspherical: no relator can
be a proper power. With this convention, asphericity of a presentation
implies asphericity of the Cayley $2$-complex (\cite{Gersten87}, p.~32),
hence (by Hurewicz' Theorem) cohomological dimension at most $2$ and
hence (\cite{B}, p.~187) torsion-freeness.

\bigskip

\begin{dem}[  of the theorem]
Let $\Gamma$ be a reduced labelled graph. The group under consideration
is defined by the presentation
$\presgroup{S}{R}$ where $R$ is the set of all words read along cycles of
$\Gamma$. However, taking all words is not necessary: the group presented by
$\presgroup{S}{R}$ will be the same if we take not all cycles but only
a generating set of cycles.

The fundamental group of the graph $\Gamma$ is a free group. Let
$\mathcal{C}$ be a finite generating set of $\pi_1(\Gamma)$ (maybe not
standard). Let $R$ be the set of words read on the cycles in
$\mathcal{C}$.

Add $2$-faces to $\Gamma$ in the following way: for each cycle in
$\mathcal{C}$, glue a disk bordering this cycle.  Denote by $\Gamma_2$
this $2$-complex; it depends on the choice of $\mathcal{C}$, or
equivalently on $R$. 

As the cycles in $\mathcal{C}$ generate all cycles, $\Gamma_2$ is simply
connected.  Note that if $\mathcal{C}$ happens to be taken standard, as
will sometimes be the case below, then $\Gamma_2$ has no homotopy in
degree $2$.

By our definitions above (Definition~\ref{def:tileofcomplex}), a tile of $\Gamma_2$ is a topological disk
whose \tempgras{boundary} is labelled by some word of $R$.

\bigskip

We are going to show that there exists a constant $C>0$ such that any
simply connected van Kampen-reduced puzzle $D$ with respect to the tiles
of $\Gamma_2$ satisfies a linear isoperimetric inequality $\abs{\d D}\geq
C\abs{D}$ where $\abs{\d D}$ is the \tempgras{boundary} length of $D$ and
$\abs{D}$ is the number of faces of $D$.  This implies hyperbolicity (see
for example~\cite{S}).

We can safely assume that all edges of $D$ lie on the \tempgras{contour} of some
face (roughly speaking, there are no ``filaments''). Indeed, filaments
only improve isoperimetry.
Generally speaking, in what follows we will never mention the possible
occurrence of filaments, their treatment being immediate.

\begin{rem}
\label{rem:uniquelift}
The $1/6$ assumption on pieces implies that no two distinct cycles
of $\Gamma$ bear the same word.
\end{rem}

Let $e$ be an internal\footnote{i.e.\ not on the \tempgras{boundary}} edge of
$D$, adjacent\footnote{We say that two faces $f_1, f_2$ of a $2$-complex are \emph{adjacent}
along edge $e$ (or simply \emph{adjacent} if the mention of $e$ is
unnecessary) if either $f_1\neq f_2$ and $e$ belongs to the
\tempgras{contour} of both $f_1$ and $f_2$, or $f_1=f_2$ and $e$ is included twice
in the \tempgras{contour} of $f_1$.} to
faces $f_1$ and $f_2$. 
As $D$ is a puzzle over the tiles of $\Gamma_2$,
there are faces $f'_1$ and $f'_2$ of $\Gamma_2$ bearing the same \tempgras{contour}
words as $f_1$ and $f_2$ respectively (maybe up to inversion).  These
faces are unique by Remark~\ref{rem:uniquelift}.

The edge $e$ belongs to the \tempgras{contour} of both $f_1$ and $f_2$
and thus can be lifted in $\Gamma_2$ either in $f'_1$ or in $f'_2$. Say
$e$ is an \emph{edge originating from $\Gamma_2$} if these two lifts
coincide, so that in $\Gamma_2$, the two faces at play are adjacent along
the same edge as they are in $D$.

Any labelled complex with respect to the tiles of $\Gamma_2$, all
internal edges of which originate from $\Gamma_2$, can thus be lifted to
$\Gamma_2$ by lifting each of its edges. This lifting is unique by
Remark~\ref{rem:uniquelift}.

Note that $D$ is van Kampen-reduced if and only if there is no edge $e$
originating from $\Gamma_2$ and adjacent to faces $f_1$, $f_2$ such that
$f'_1=f'_2$.

\bigskip

We work by first proving the isoperimetric inequality for puzzles having all
edges originating from $\Gamma_2$. Second, we will decompose the puzzle $D$
into ``parts'' having all their edges originating from $\Gamma_2$ and show that
these parts are in $1/6$ small cancellation with each other. Then we will
use ordinary small cancellation theory to conclude.

We begin by proving what we want for some particular choice of $R$.

\begin{lem}\label{lemmain}
Let $\Delta=\diam(\Gamma)$. Suppose that $\mathcal{C}$ was chosen to be the set of 
closed paths embedded (or immersed) in $\Gamma$ of length at most $3\Delta$. Then, for
any closed path in $\Gamma$ labelling a reduced word $w$, there exists a
simply connected puzzle with \tempgras{boundary} word $w$, with tiles having their
\tempgras{boundary} words in $R$, all edges of which originate from $\Gamma_2$, and
with at most $3\abs{w}/g$ tiles.
\end{lem}

\begin{dem}[ of Lemma~\ref{lemmain}]
If $\abs{w} \leq 2\Delta$ then by definition of $R$ there exists a one-tile puzzle
spanning $w$, and as $\abs{w}\geq g$ the conclusion holds.
Show by induction on $n$ that if $\abs{w}\leq n\Delta$ there
exists a simply connected puzzle $D$ spanning $w$ with at most $n$ tiles. This is true
for $n=2$. Suppose this is true up to $n\Delta$ and suppose that
$2\Delta\leq \abs{w}\leq (n+1)\Delta$.

Let $w=w'w''$ where $\abs{w'}=2\Delta$. As the diameter of $\Gamma$ is
$\Delta$, there exists a path in $\Gamma$ labelling a word $x$ joining
the endpoints of $w'$, with $\abs{x}\leq \Delta$. So $w'x^{-1}$ is read
on a cycle of $\Gamma$ of length at most $3\Delta$, hence (its reduction)
belongs to $R$. Now $xw''$ is a word read on a cycle of $\Gamma$, of
length at most $\abs{w}-\Delta\leq n\Delta$. So there is a puzzle with at
most $n$ tiles spanning $xw''$. Gluing this puzzle with the tile spanning
$w'x^{-1}$ along the $x$-sides provides the desired puzzle. (Note that this gluing
occurs in $\Gamma_2$, so that edges of the resulting puzzle originate from
$\Gamma_2$.)

So for any $w$ we can find a puzzle spanning it with at most
$1+\abs{w}/\Delta$ tiles. As $\Delta\geq g/2$ and as $\abs{w}\geq g$, we
have $1+\abs{w}/\Delta\leq 1+2\abs{w}/g \leq 3\abs{w}/g$.
\end{dem}

\begin{cor}\label{lemmult}
For any choice of $\mathcal{C}$, there exists a constant $\alpha>0$ such that any
minimal simply connected puzzle $D$ with
respect to the tiles of $\Gamma_2$ all internal edges of which originate
from
$\Gamma_2$ satisfies the isoperimetric inequality $\abs{\d D}\geq
\alpha\abs{D}$.
%
\end{cor}

\begin{dem}[ of Corollary~\ref{lemmult}]
Indeed, the existence of an isoperimetric constant for minimal diagrams
does not depend on the finite presentation, hence the result when
$\mathcal{C}$ is finite. This also holds for infinite $\mathcal{C}$
since any infinite family of cycles in the finite graph $\Gamma$ contains a finite generating subfamily.
%
\end{dem}

These last affirmations only express in terms of diagrams the fact that
the fundamental group of $\Gamma$, which is free hence hyperbolic, is
generated by the cycles of $\Gamma$ of length at most $3\Delta$ (w.r.t.\
some basepoint).

\bigskip

The next lemma is just ordinary small cancellation theory (see for
example the appendix of~\cite{GH}, or~\cite{LS}), stated in the form
we need.
Note that usually, the definition of small cancellation involves pieces
of relative size less than $\lambda$ with $\lambda\leq 1/6$. Here we use
pieces of relative size at most $\lambda$ with $\lambda<1/6$. This is
less well-suited for treatment of infinite presentations (which we do not
consider) but allows lighter notation for the isoperimetric constant
$1-6\lambda>0$ and the Greendlinger constant $1-3\lambda>1/2$.

\begin{lem}\label{lemsc}
Let $R$ be a set of simply connected reduced tiles. Suppose that any piece with
respect to two tiles $t,t'\in R$ is a word of length at most $\lambda$ times 
the smallest \tempgras{boundary} length of $t$ and $t'$, for some constant
$\lambda<1/6$.

Then any simply connected van Kampen-reduced puzzle $D$ with respect to the tiles of $R$
satisfies the following properties.
\begin{enumerate}
\item If $D$ has at least two faces, the reduction $w$ of the \tempgras{boundary} word of $D$
contains two disjoint subwords $w_1$, $w_2$, with $w_1$ (resp.\ $w_2$)
subword of the \tempgras{boundary} word of some tile $t_1$ (resp.\ $t_2$) of $D$,
with length at least $(1-3\lambda)>\frac12$ times the \tempgras{boundary} length of
$t_1$ (resp.\ $t_2$).
\item
The word $w$ is not a
proper subword of the \tempgras{boundary} word of some tile.
\item
The \tempgras{boundary} length $\abs{\d D}$ is at
least $1-6\lambda$ times the sum of the lengths of the faces of $D$, and
at least the \tempgras{boundary} length of the largest tile it contains.
\end{enumerate}

Moreover, there is no spherical van Kampen-reduced puzzle with respect to
these tiles.
\end{lem}

\begin{cor}\label{nohole}
Let $R$ be a set of (not necessarily simply connected) reduced tiles. Suppose
that any piece with respect to two tiles $t,t'\in R$ is a word of
length at most $\lambda$ times the smallest length of the \tempgras{boundary}
component of $t$ and $t'$ it immerses in, for some constant $\lambda<1/6$.

Then, any simply connected puzzle with respect to this set of tiles
contains only simply connected tiles.
\end{cor}

\begin{dem}[ of the corollary]

Let $D$ be a simply connected puzzle with respect to $R$. Let $t$ be a
non-simply connected tile in $D$. We can suppose that $t$ is deepest,
that is, that the bounded components of the complement of $t$ contain no
other non-simply connected tile.

The interior of $t$ is embedded in the plane and is homeomorphic to a
disk with some finite number $n$ of holes. Since $D$ is simply connected
any such hole is filled with a subpuzzle. So let $D'_1,\ldots,D'_n$ be
the subpuzzles filling the bounded connected components of the complement
of the interior of $t$.  Each $D'_i$ is simply connected, since the
bounded connected components of the complement of a connected set in the
plane are simply connected.  Let us work with $D'_1$. In case $D'_1$ is
not van Kampen-reduced we replace it by its van Kampen-reduction (which
does not change its boundary word, so it can still be glued to one of the
holes of $t$).

The \tempgras{boundary} of $D'_1$ may not be embedded in the plane.
However, it is immersed, since the word read on it is the word read on
one of the interior boundaries of $t$, and this word is reduced.

The component $D'_1$ is a connected simply connected puzzle. Its image in
the plane is the union of closed sets $D''_1,\ldots,D''_q$ such that each
$D''_i$ is either a topological closed disk or a topological closed
segment (``filament''), and the $D''_i$'s intersect at a finite number of
points.
By construction,
each $D''_i$ which is a disk is a puzzle.

\begin{center}
\includegraphics{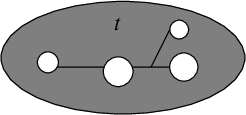}
\end{center}

Suppose that $D''_i$ is a segment. Then each of its endpoints belongs to
some $D''_j$ with $j\neq i$. Indeed, otherwise the \tempgras{boundary} of $D'_1$
would not be immersed.

Construct a graph $T$ embedded in the plane in the following way.  For
each $D''_i$ which is a disk, define a family of segments $T_i$ as
follows: Choose a point $p_0$ in the interior of $D''_i$. There are a
finite number of points $p_1,\ldots,p_r$ on the \tempgras{boundary} of
$D''_i$ such that $p_j$ belongs to some $D''_k$ for $k\neq i$. Now define
$T_i$ to be made of the union of segments $p_0p_j\subset D''_i$ for
$1\leq j \leq r$.  Now define $T$ to be the union of all $D''_i$ for
those $1\leq i \leq q$ for which $D''_i$ is a segment, plus the union of
all $T_i$'s for those $1\leq i \leq q$ for which $D''_i$ is a disk.

By construction, $T$ is connected since $D'_1$ is.

For each $i$ such that $D''_i$ is a disk, $D''_i$ retracts onto $T_i$
preserving the points $p_1,\ldots,p_r$. So $D'_1$ retracts onto $T$, and
in particular $T$ is simply connected since $D'_1$ is. So $T$ is a tree.
It is non-empty since $D'_1$ is (but maybe reduced to a point if $D'_1$
is a topological disk).

Now consider some leaf of $T$. Since any endpoint of any $D''_i$ which is
a segment belongs to some $D''_j$ with $j\neq i$ (since $\partial D'_1$
is immersed as we saw above), a leaf of $T$ cannot belong to a $D''_i$
which is a segment. So a leaf of $T$ belongs to some $T_i$ constructed
from some $D''_i$ which is a disk. By definition of $T_i$, this means
that $D''_i$ intersects with at most one other $D''_j$ with $j\neq i$.

Now $D''_i$ is a puzzle which is a topological disk. As we supposed that
$t$ was taken a deepest non-simply connected tile, $D''_i$ contains only
simply connected tiles. So we can apply Lemma~\ref{lemsc}:  there exist
two tiles $t'$, $t''$ in $D''_i$ and two subwords $w'$, $w''$ of the
\tempgras{boundary} word of $D''_i$ such that $w'$ (resp.\ $w''$) is a subword of
the \tempgras{boundary} word of $t'$ (resp.\ $t''$) of length at least one half the
\tempgras{boundary} length of $t'$ (resp.\ $t''$). As $D''_i$ has at most one point
of intersection with the other $D''_j$ for $j \neq i$, at least one of
$w'$ and $w''$ is a subword of the \tempgras{boundary} of $D'_1$. But a \tempgras{boundary} word of
$D'_1$ is a \tempgras{boundary} word of the tile $t$, and so $t$ shares with $t'$ or
$t''$ a word of length at least one half the \tempgras{boundary} length of $t'$ or
$t''$, which contradicts the small cancellation assumption.
\end{dem}

Back to our simply connected van Kampen-reduced minimal puzzle $D$ with
tiles in $\Gamma_2$. A puzzle is built by taking the disjoint union of
all its tiles and gluing them along the internal edges.

First, define a disjoint union of puzzles $D'$ by taking the disjoint
union of all tiles of $D$ and gluing them along the internal edges of
$D$ originating from $\Gamma_2$. All internal edges of $D'$ originate from
$\Gamma_2$.

As $D$ is van Kampen-reduced, $D'$ is as well.

Let $D_i$, $i=1,\ldots,n$ be the connected components of $D'$. They form
a partition of $D$. The puzzle $D$ is obtained by gluing these components
along the internal edges of $D$ not originating from $\Gamma_2$.

It may be the case that the \tempgras{boundary} word of some $D_i$ is not reduced.
This means that there is a vertex on the \tempgras{boundary} of $D_i$ which is the
origin of two (oriented) edges bearing the same vertex. We will modify
$D$ in order to avoid this. Suppose some $D_i$ has non-reduced
\tempgras{boundary} word
and consider two edges $e_1, e_2$ of $D$ responsible for this: $e_1$ and
$e_2$ are two consecutive edges with inverse labels. These edges
are either \tempgras{boundary} edges of $D$ or internal edges. In the latter case
this means that $D_i$ is to be glued to some $D_j$. We treat only this
latter case as the other one is even simpler.

Make the following transformation of $D$: do not glue any more edge $e_1$
of $D_i$ with edge $e_1$ of $D_j$, neither edge $e_2$
of $D_i$ with edge $e_2$ of $D_j$, but rather glue edges $e_1$ and $e_2$
of $D_i$, as well as edges $e_1$ and $e_2$ of $D_j$, as in the
following picture. This is possible since by definition $e_1$ and $e_2$
bear inverse labels.

\begin{center}
\includegraphics[scale=0.7]{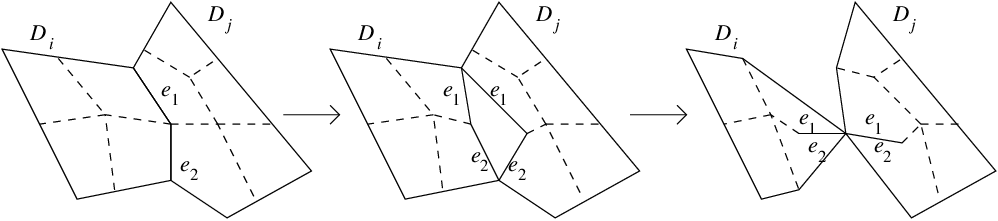}
\end{center}

This kind of operation has been studied and termed \emph{diamond move} in
\cite{CH}. The case when the central point has valency greater than $2$
(i.e.\ when more than two $D_i$'s meet at this point) is treated
similarly.

Since $\Gamma_2$ is reduced, the lifts to $\Gamma_2$ of the edges $e_1$ and
$e_2$ of $D_i$ are the same edge of $\Gamma_2$. This shows that the
transformation above preserves the fact that all edges of $D_i$ and of
$D_j$ originate from $\Gamma_2$.

The resulting puzzle (denoted $D$ again) has the same number of faces as
before, and no more \tempgras{boundary} edges. Thus, proving isoperimetry for the
modified puzzle will imply isoperimetry for the original one as well.
So we can safely assume that the \tempgras{boundary} words of the $D_i$'s are
reduced.

Now consider $D$ as a puzzle with the $D_i$'s as tiles. (More precisely,
if we erase from $D$ all internal edges originating from $\Gamma_2$ then
we obtain a puzzle each tile of which is the tile associated to the
one-face complex obtained from some $D_i$ by erasing all internal edges
originating from $\Gamma_2$.) This is a van Kampen-reduced puzzle, since
if $D_i$ and $D_j$ are in reduction position this means that they lift to
the same subcomplex of $\Gamma_2$ and share an edge originating from
$\Gamma_2$, which contradicts their definition.  Note that these tiles
are not necessarily simply connected. 

These tiles satisfy the condition of Corollary~\ref{nohole}. Indeed,
suppose that two tiles $D_i$, $D_j$ (with maybe $i=j$ in which case two
parts of the \tempgras{boundary} of the same tile are glued) are to be
glued along a common (reduced!) word $w$.  By definition of the $D_i$'s,
the edges making up $w$ do not originate from $\Gamma_2$.

As the edges of $D_i$ originate from $\Gamma_2$, there is a lift
$\phi_i:D_i\rightarrow \Gamma_2$ (as noted above).
Consider the two lifts $\phi_i(w)$ and $\phi_j(w)$. As the edges making
up $w$ do not originate from $\Gamma_2$, these two lifts are different. As
$w$ is reduced these lifts are immersions.  So $w$ is a piece. By
assumption the length of $w$ is at most $\Lambda<g/6$.

Now as $D_i$ lifts to $\Gamma_2$, any \tempgras{boundary} component of $D_i$ goes to a
closed path in $\Gamma$. This proves that the length
of any \tempgras{boundary} component of $D_i$ is at least $g$.

So the tiles $D_i$ satisfy the small cancellation condition with
$\lambda=\Lambda/g<1/6$. As they are tiles of a simply connected puzzle,
by Corollary~\ref{nohole} they are simply connected.

Then by Lemma~\ref{lemsc}, the \tempgras{boundary} of $D$ is at
least $1-6\lambda$ times the sum of the \tempgras{boundary} lengths of the $D_i$'s
(
considered as tiles).
Since $D$ is minimal, each $D_i$ is as well, and as $D_i$ is simply
connected, by Corollary~\ref{lemmult} it satisfies the isoperimetric inequality
$\abs{\d D_i}\geq \alpha \abs{D_i}$. So
\[
\abs{\d D}\geq (1-6\lambda) \sum \abs{\d D_i}\geq \alpha (1-6\lambda)
\sum \abs{D_i}=\alpha (1-6\lambda)\abs{D}
\]
which shows the isoperimetric inequality for $D$, hence hyperbolicity.

\bigskip

For asphericity and the cohomological dimension (hence torsion-freeness),
suppose that $\mathcal{C}$ is standard (so that $\Gamma_2$ is
contractible) and that there exists a van Kampen-reduced spherical puzzle
$D$, which we can assume to be inclusion-minimal in the sense that it
contains no spherical subpuzzle. Define the $D_i$'s as above.  Either
some $D_i$ is spherical, in which case $D=D_i$ by inclusion-minimality of
$D$, or all $D_i$'s have non-empty boundary words. The former is ruled out
by the following lemma:

\begin{lem}
Suppose that the set of paths read along faces of $\Gamma_2$ is standard.
Let $D$ be a non-empty spherical puzzle all edges of which originate from $\Gamma_2$.
Then $D$ is not van-Kampen reduced.
\end{lem}

\begin{dem}[ of the lemma]
Let $T$ be a maximal tree of $\Gamma$ witnessing for
standardness of the family of cycles. Homotope $T$
to a point. This turns $\Gamma_2$ into a bouquet of circles with a face in
each circle. Similarly, homotope to a point any edge of $D$ coming from a
suppressed edge of $\Gamma$. This way we turn $D$ into a spherical van
Kampen diagram with respect to the presentation of the fundamental group
of $\Gamma_2$ (i.e.\ the trivial group) by 
$\presgroup{c_1,\ldots,c_n}{c_1=e,\ldots,c_n=e}$. But there is no reduced
spherical van Kampen diagram with respect to this presentation, as can
immediately be checked.
\end{dem}

Since by definition each $D_i$ lifts to $\Gamma_2$ and since $D$ (hence
each $D_i$) is van Kampen-reduced, the lemma implies that no $D_i$ is
spherical. Hence the $D_i$'s have non-trivial boundary words. So $D$ can
be viewed as a spherical puzzle with the boundary words of the $D_i$'s as tiles.
But we saw above that the $D_i$'s (viewed as tiles) satisfy the small
cancellation condition. So by Lemma~\ref{lemsc} there is no spherical
van Kampen-reduced puzzle w.r.t.\ these tiles.

The computation of the Euler characteristic immediately follows, using
that the cohomological dimension is at most $2$.

The last assertions of the theorem follow easily from the assertions
of Lemma~\ref{lemsc}. The smallest relation in the group presented by
$\presgroup{S}{R}$ is the \tempgras{boundary} length of the smallest
non-trivial puzzle, which by Lemma~\ref{lemsc} is at least the smallest
\tempgras{boundary} length of the $D_i$'s, which is at least the girth $g$.
Similarly, any reduced word representing the trivial element in the group
is read on the \tempgras{boundary} of a van-Kampen reduced simply
connected puzzle, thus
contains as a subword at least one half of the \tempgras{boundary} word
of some $D_i$.

\bigskip

{\footnotesize \textit{(Note:
The version of this text published in Bull.\ Belg.\ Math.\ Soc.\ contains
a mistake in this part of the argument, as it used that $x'$ was
geodesic. This was pointed to me by Mikhail Ostrovskii. Below is a corrected
version.)
}}

For the isometric embedding of $\Gamma$ in the Cayley graph of the
group, suppose that some geodesic path $p$ in the graph (or in
$\Gamma_2$) labelling a word $x$
is equal to a shorter word $y$ in the quotient. This means that there
exists a van Kampen-reduced puzzle $D$ with \tempgras{boundary} word $xy^{-1}$, made up of tiles with
cycles of $\Gamma$ as \tempgras{boundary} words.

Now $x$ is the word read on a path $p_1$ in the \tempgras{boundary} word of $D$,
which lifts to the geodesic path $p$ in $\Gamma_2$ labelling $x$ as well.
Let $f$ be a face of $D$ which intersects $p_1$ along at least one edge. We
say that $f$ \emph{originates from $\Gamma_2$ together with $x$} if the
lift from $f$ to $\Gamma_2$ coincides with the lift $p_1\rightarrow
p$ on the intersection of $f$ with $p_1$.

We are going to recursively remove all faces of $D$ which originate from
$\Gamma_2$ together with $x$, as follows. Let $f$ be such a face of $D$,
and assume it shares an edge $e$ along with $x$, so that $x=x_1ex_2$ and
the boundary of $D$ is $ew$. Define the path $x'=x_1w^{-1}x_2$ in the
diagram $D$, and remove face $f$ from the diagram $D$. This defines a new
diagram $D'$ with boundary $x'y^{-1}$. Note that $x'$ may not be reduced.

By construction, $x'$ lifts to $\Gamma_2$
together with $x$, and their lifts to $\Gamma_2$ have the same endpoints.
In particular, $\abs{x'}\geq \abs{x}$ since $x$ is geodesic in $\Gamma$.

Repeat this process until there are no faces of $D$ that originate from
$\Gamma_2$ together with $x'$. At the end of this process, we still have
that $x'$ lifts to $\Gamma_2$ together with $x$, and $\abs{x'}\geq
\abs{x}$. 
In
the following picture in which black cells represent tiles originating
from $\Gamma_2$ together with $x$.

\begin{center}
\includegraphics[scale=0.5]{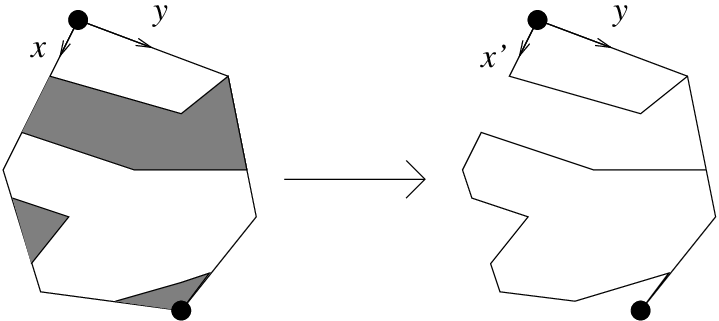}
\end{center}

At this point, if $x'=y$ (there are no faces left in the diagram), we are
done, because $\abs{y}=\abs{x'}\geq \abs{x}$ as needed.

If some faces are left, we have a (possibly non-reduced) puzzle
with boundary $x'y^{-1}$.
By construction, no faces of $D'$ lift to $\Gamma_2$ together
with $x'$. This means that the intersection of $x'$ with any face of $D'$
is a \emph{piece} in $\Gamma$.

By assumption, $y$ was reduced, but
there may be cancellations within $x'$ or between $x'$ and $y$. Reduce
$D'$, first by removing any filaments in $x'$, then by ``folding
in'' any inverse consecutive edges of $x'y^{-1}$ as in the ``diamond
moves'' above. Let $w_1$ and $w_2$ be the common initial and final
segments between $x'$ and $y$ (if any); after reduction, we are left with
a puzzle $D''$ with boundary word $x''(y')^{-1}$, where $x'=w_1x''w_2$,
$y=w_1y'w_2$, and $x''(y')^{-1}$ is cyclically reduced.

Now (if $\Gamma$ contains no filaments) $x''$ is part of some cycle of
$\Gamma$
labelled by $x''z$.
%
As $\Gamma_2$ is simply connected, there is a van Kampen-reduced puzzle $D_2$ with boundary
word $x''z$ and which globally lifts to $\Gamma_2$ (all its edges originate
from $\Gamma_2$).

Define a new puzzle $D'''$ by gluing $D_2$ and $D''$ along
the word $x''$. This is a puzzle bordering $zy'$. It is van
Kampen-reduced since $D_2$ and $D''$ are van Kampen-reduced and
since there is no cancellation between $D_2$ and $D''$ (otherwise there
would be a tile of $D''$ originating from $\Gamma_2$ together with $x''$).

Now consider, as above, the partition $D'''=\cup D'''_i$
where the
$D'''_i$ are maximal parts lifting to $\Gamma_2$. Since no tile of $D''$
adjacent to $x''$ originates from $\Gamma_2$ together with $x''$, $D_2$ is
exactly one of the $D'''_i$.

By Lemma~\ref{lemsc}, the \tempgras{boundary} length
$\abs{z}+\abs{y'}$ of $D'''$ is at least the \tempgras{boundary} length of
any $D'''_i$. In particular, it is at least the \tempgras{boundary}
length of $D_2$, which is $\abs{z}+\abs{x''}$. This proves that
$\abs{z}+\abs{y'}\geq \abs{z}+\abs{x''}$, and therefore,
$\abs{y'}\geq\abs{x''}$ so that $\abs{y}\geq \abs{x'}\geq \abs{x}$, as needed.

\bigskip

This proves the theorem.
\end{dem}

\section{Further remarks}

\begin{rem}
The proof above gives an explicit isoperimetric constant when the set of
relators taken is the set of all words read on cycles of the graph of
length at most three times the diameter: in this case, any
minimal simply connected puzzle satisfies the
isoperimetric inequality
\[
\abs{\d D}\geq g(1-6\Lambda/g)\abs{D}/3
\]
This explicit isoperimetric constant growing linearly with $g$ (i.e.~``homogeneous'') can be
very useful if one wants to apply such theorems as the local-global
hyperbolic principle, which requires the isoperimetric constant to grow
linearly with the sizes of the relators.
\end{rem}

\begin{rem}
The assumption that $\Gamma$ is reduced can be relaxed a little bit,
provided that some quasi-geodesicity assumption is granted, and that the
definition of a piece is emended.

Redefine a \emph{piece} to be a couple of words $(w_1,w_2)$ such that
both immerse in $\Gamma$ and such that $w_1=w_2$ in the free group. The
\emph{length} of a piece $(w_1,w_2)$ is the maximal length of $w_1$ and
$w_2$.

There are trivial pieces, for example if $w_1=w_2$ and both have the
same immersion. However, forbidding this is not enough: for example, if a
word of the form $aa^{-1}w$ immerses in the graph, then $(aa^{-1}w, w)$
will be a piece.

A \emph{trivial piece} is a piece $(w_1,w_2)$ such that
there exists a path $p$ in $\Gamma$ joining the beginning of the
immersion of $w_1$ to the beginning of the immersion of $w_2$ such that
$p$ is labelled with a word equal to $e$ in the free group.

The new theorem is as follows.
\end{rem}

\begin{thm}[ (M.~Gromov)]
\label{thm:variantmain}
Let $\Gamma$ be a finite non-filamenteous labelled graph. Let $R$ be the set of words read on
all cycles of $\Gamma$ (or on a generating family of cycles). Let $g$ be the girth of $\Gamma$ and $\Lambda$ be the
length of the longest non-trivial piece of $\Gamma$.

Suppose that $\lambda=\Lambda/g$ is less than $1/6$.

Suppose that there exist a constant $A>0$ such that any word $w$ immersed
in $\Gamma$ of length at least $L$ satisfies $\norm{w}\geq
A(\abs{w}-L)$ for some $L<(1-6\lambda)g/2$.

Then the presentation $\presgroup{S}{R}$ defines a hyperbolic,
infinite, torsion-free group $G$, and (if $R$ arises from a standard
family of cycles) this presentation is aspherical
(hence the cohomological dimension of $G$ is at most $2$). Moreover, the
natural map of labelled graphs from $\Gamma$ to the Cayley graph of $G$
is a $(1/A,AL)$-quasi-isometry. The shortest relation of
$G$ is of length at least $Ag/2$, and any reduced word equal to $e$ in $G$
contains as a subword the reduction of at least one half of a word read
on a cycle of $\Gamma$.
\end{thm}

(Here $\norm{w}$ is the length of the reduction of $w$; besides,
in accordance with~\cite{GH}, by a $(\lambda,c)$-quasi-isometry we wean
a map $f$ such that $d(x,y)/\lambda-c\leq d(f(x),f(y))\leq \lambda
d(x,y)+c$.)

\begin{rem}
The same kind of theorem holds if we use the $C(7)$ condition instead of the
$C'(1/6)$ condition, but in this case there is no control on the radius
of injectivity (shortest relation length).
\end{rem}

\begin{rem}
Using the techniques in~\cite{D} or~\cite{O}, the same kind of theorem should
hold starting with any torsion-free hyperbolic group instead of the free
group, provided that the
girth of the graph is large enough w.r.t.\ the hyperbolicity constant,
and that the labelling is quasi-geodesic. See~\cite{Ollexp}.
\end{rem}

\begin{rem}
Theorem~\ref{thm:main} can be extended when the graph is infinite, in which
case we get a direct limit of torsion-free, dimension-$2$ hyperbolic
groups (but generally not
hyperbolic), in which the conclusions of small cancellation theory still
hold but with the isoperimetric constant for van Kampen diagrams tending
to $0$. In this case the small cancellation assumption reads: any piece
has length less than $1/6$ times the minimal length of a cycle on which
it appears.
\end{rem}

\end{document}